\documentclass[12pt,a4paper]{article}
\usepackage{amsmath,amsfonts,amssymb,amscd}

\usepackage[a4paper]{geometry}

\usepackage{ulem}
\usepackage{xcolor}
\usepackage{graphics}

\usepackage{hyperref}

\def\gr{\operatorname{gr}}

\def\red{\operatorname{red}}

\def\GL{\operatorname{GL}}

 \newtheorem{thm}{Theorem}[section]
 \newtheorem{cor}[thm]{Corollary}
 \newtheorem{lem}[thm]{Lemma}
 \newtheorem{prop}[thm]{Proposition}

\newcounter{th}

\newcounter{prop}

\newcounter{lem}

\newcounter{de}
\def\de{\refstepcounter{de}{\bf \noindent{Definition} \arabic{de}. }}

\newcounter{ex}

\begin{document}

\begin{center}
{\Large On complex-analytic $1|3$-dimensional supermanifolds associated with $\mathbb{CP}^1$.}
\footnote[1]{Supported by Fonds National de la Recherche Luxembourg\\

}

\bigskip

{\bf E.G.~Vishnyakova}\\
\end{center}

\bigskip

\begin{abstract} We obtain a classification  up to isomorphism of
complex-analytic supermanifolds with underlying space $\mathbb{CP}^1$
 of dimension $1|3$ with retract $(k,k,k)$, where $k\in \mathbb{Z}$. More precisely, we prove
that classes of isomorphic complex-analytic supermanifolds of dimension $1|3$ with retract $(k,k,k)$
are in one-to-one correspondence with points of the following set:
$$
\mathbf{Gr}_{4k-4,3} \cup \mathbf{Gr}_{4k-4,2} \cup \mathbf{Gr}_{4k-4,1} \cup \mathbf{Gr}_{4k-4,0}
$$
for $k \geq 2$. For $k < 2$ all such supermanifolds are isomorphic to their retract $(k,k,k)$.

\end{abstract}

\section{Introduction.} A classical result is that we can assign the holomorphic vector bundle,
  so called retract, to each complex-analytic supermanifold (see Section 2 for more details). Assume that the underlying space of a
   complex-analytic supermanifold is $\mathbb {CP}^1$. By the Birkhoff-Grothendieck Theorem
any vector bundle of rank $m$ over $\mathbb {CP}^1$ is isomorphic
  to the direct sum of $m$ line bundles: $\mathbf E \simeq \bigoplus\limits_{i=1}^m L(k_i)$,
  where $k_i\in \mathbb{Z}$. We obtain a
  classification up to isomorphism of complex-analytic supermanifolds of dimension
$1|3$ with underlying space $\mathbb {CP}^1$  and with retract $L(k)\oplus L(k)\oplus L(k)$,
where $k\in \mathbb{Z}$. In addition, we give a classification up to isomorphism of complex-analytic
supermanifolds of dimension $1|2$ with underlying space $\mathbb {CP}^1$.

 The paper is structured as follows. In Section $2$ we explain the
idea of the classification.
  In Section $3$ we do all necessary preparations.
 The classification up
to isomorphism of complex-analytic supermanifolds of
dimension $1|3$ with underlying space $\mathbb{CP}^1$ and  with retract $(k,k,k)$ is
obtained in Section $4$. The last section is devoted to the classification up to
isomorphism of complex-analytic supermanifolds of dimension $1|2$ with underlying space $\mathbb{CP}^1$.

The study of complex-analytic supermanifolds
 with underlying space $\mathbb {CP}^1$ was started
in \cite{Bunegina_Oni}. There the classification of homogeneous
complex-analytic supermanifolds of dimension $1|m$, $m\le 3$, up to
isomorphism was given. It was proven that in the case $m=2$ there
exists only one non-split homogeneous supermanifold constructed by
P.~Green \cite{Green} and V.P.~Palamodov \cite{Ber}. For $m = 3$ it
was shown that there exists a series of non-split homogeneous
supermanifolds, parameterized by $k=0,2,3,\cdots$.

In \cite{Vish} we studied even-homogeneous supermanifold, i.e. supermanifolds which possess transitive actions of Lie
groups. It was shown that there exists a series of non-split even-homo\-ge\-ne\-ous
supermanifolds, parameterized by elements in $\mathbb{Z}\times
\mathbb{Z}$, three series of non-split even-homogeneous
supermanifolds, parameterized by elements of $\mathbb{Z}$, and a
finite set of exceptional supermanifolds.

\section{Classification of supermanifolds, main definitions}

We will use the word "supermanifold" \,in the sense of Berezin --
Leites \cite{BL, ley}, see also \cite{Bunegina_Oni}. All the time, we will be interested in
the complex-analytic version of the theory. We begin with main definitions.

Recall that a {\it complex superdomain of dimension $n|m$} is a $\mathbb{Z}_2$-graded
ringed space of the form $(U, \mathcal{F}_U \otimes \bigwedge(m) )$,
where $\mathcal{F}_U$ is the sheaf of holomorphic functions on an open set $U\subset \mathbb{C}^n$ and
$ \bigwedge(m)$ is the exterior (or Grassmann) algebra with $m$ generators.
\medskip

\de A {\it complex-analytic supermanifold} of dimension $n|m$ is a $\mathbb{Z}_2$-graded locally ringed space that
is locally isomorphic to a complex superdomain of dimension $n|m$.

\medskip

Let $\mathcal{M} = (\mathcal{M}_0,{\mathcal
O}_{\mathcal{M}})$ be a supermanifold and
$$
\mathcal{J}_{\mathcal{M}} = ({\mathcal
O}_{\mathcal{M}})_{\bar 1} + (({\mathcal
O}_{\mathcal{M}})_{\bar 1})^2
$$
be the subsheaf of ideals generated
by odd
elements in ${\mathcal O}_{\mathcal{M}}$. We put $\mathcal{F}_{\mathcal{M}}:= {\mathcal
O}_{\mathcal{M}}/\mathcal{J}_{\mathcal{M}}$. Then $(\mathcal{M}_0,
\mathcal{F}_{\mathcal{M}})$ is a usual complex-analytic manifold, it
is called the {\it reduction} or {\it underlying space} of $\mathcal{M}$. Usually we will
write $\mathcal{M}_0$ instead of $(\mathcal{M}_0,
\mathcal{F}_{\mathcal{M}})$.  Denote by ${\mathcal T}_{\mathcal{M}}$
the {\it tangent sheaf} or the {\it sheaf of vector fields} of
$\mathcal{M}$. In other words, ${\mathcal T}_{\mathcal{M}}$ is the
sheaf of derivations of the structure sheaf ${\mathcal
O}_{\mathcal{M}}$. Since the sheaf ${\mathcal O}_{\mathcal{M}}$ is
$\mathbb{Z}_2$-graded, the tangent sheaf ${\mathcal
T}_{\mathcal{M}}$ also possesses the induced $\mathbb{Z}_2$-grading,
i.e. there is the natural decomposition ${\mathcal T}_{\mathcal{M}}
= ({\mathcal T}_{\mathcal{M}})_{\bar 0} \oplus ({\mathcal
T}_{\mathcal{M}})_{\bar 1}$.

Let $\mathcal{M}_0$ be a complex-analytic manifold and let $\mathbf{E}$ be a holomorphic vector bundle over $\mathcal{M}_0$. Denote by
$\mathcal{E}$ the sheaf of holomorphic
sections of $\mathbf{E}$. Then the ringed space
$(\mathcal{M}_0,\bigwedge\mathcal{E})$ is a supermanifold. In this case $\dim\,\mathcal{M} = n|m$, where $n =
\dim \mathcal{M}_0$ and $m$ is the rank of ${\mathbf E}$.

\medskip
\de
 A supermanifold $(\mathcal{M}_0,{\mathcal
O}_{\mathcal{M}})$ is called {\it split} if
${\mathcal
O}_{\mathcal{M}}\simeq \bigwedge\mathcal{E}$ (as supermanifolds) for a
 locally free sheaf $\mathcal{E}$.
 \medskip

  It
is known that any real (smooth or real analytic) supermanifold is
split.
  The
structure sheaf ${\mathcal O}_{\mathcal{M}}$ of a split
supermanifold possesses a ${\mathbb Z}$-grading, since ${\mathcal O}_{\mathcal{M}} \simeq\bigwedge\mathcal{E}$ and the sheaf
$\bigwedge\mathcal{E}  = \bigoplus\limits_p \bigwedge^p\mathcal{E}$
is naturally ${\mathbb Z}$-graded. In other words, we have the decomposition
${\mathcal O}_{\mathcal{M}} = \bigoplus\limits_p ({\mathcal O}_{\mathcal{M}})_p$.  This ${\mathbb Z}$-grading
induces the ${\mathbb Z}$-grading in ${\mathcal T}_{\mathcal{M}}$ in the following way:
\begin{equation}\label{filtr of T}
({\mathcal T}_{\mathcal{M}})_p: = \{ v\in {\mathcal T}_{\mathcal{M}}\,  |\,  v(({\mathcal O}_{\mathcal{M}})_q)
\subset ({\mathcal O}_{\mathcal{M}})_{p+q} \,\, \text{for all} \,\,q\in \mathbb{Z}\} .
\end{equation}
In other words, we have the decomposition: ${\mathcal T}_{\mathcal{M}} =\bigoplus\limits_{p=-1}^m({\mathcal T}_{\mathcal{M}})_p$. Therefore the
superspace $H^0(\mathcal{M}_0,{\mathcal T}_{\mathcal{M}})$ is also
${\mathbb Z}$-graded. Consider the subspace
$$
{\mathrm
{End}}\,{\mathbf E}\subset H^0(\mathcal{M}_0,({\mathcal
T}_{\mathcal{M}})_0).
$$
It consists of all endo\-mor\-phisms of the vector bundle ${\mathbf
E}$ inducing the identity morphism on $\mathcal{M}_0$. Denote by
${\mathrm {Aut}}\,{\mathbf E}\subset{\mathrm {End}}\,{\mathbf E}$
the group of automorphisms of ${\mathbf E}$, i.e. the group of all
invertible endomorphisms of ${\mathbf E}$. We have the action ${\mathrm{Int}}$ of
${\mathrm {Aut}}\,{\mathbf E}$ on ${\mathcal T}_{\mathcal{M}}$ defined by
$$
{\mathrm{Int}} A: v\mapsto A v A^{-1}.
$$
Clearly, the action ${\mathrm{Int}}$ preserves
the $\mathbb{Z}$-grading (\ref{filtr of T}), therefore, we have the action of ${\mathrm
{Aut}}\,{\mathbf E}$ on $H^1(\mathcal{M}_0,{(\mathcal
T}_{\mathcal{M}})_2)$.

There is a functor denoting by $\gr$ from the category of supermanifolds
to the category of split supermanifolds. Let us describe
this construction. Let
$\mathcal{M}$ be a supermanifold and let as above
$\mathcal{J}_M\subset \mathcal{O}_{\mathcal{M}}$ be the subsheaf of ideals
generated by odd elements of $\mathcal{O}_{\mathcal{M}}$. Then by definition we have
$\gr(\mathcal{M}) = (\mathcal{M}_0,
\gr\mathcal{O}_{\mathcal{M}})$, where
$$
\gr\mathcal{O}_{\mathcal{M}}= \bigoplus_{p\geq 0} (\gr\mathcal{O}_{\mathcal{M}})_p, \quad (\gr\mathcal{O}_{\mathcal{M}})_p=
\mathcal{J}_{\mathcal{M}}^p/\mathcal{J}_{\mathcal{M}}^{p+1}, \quad
\mathcal{J}_{\mathcal{M}}^0:=\mathcal{O}_{\mathcal{M}}.
$$
In this case $(\gr\mathcal{O}_{\mathcal{M}})_1$ is a locally free sheaf and
there is a natural isomorphism of $\gr\mathcal{O}_{\mathcal{M}}$ onto $\bigwedge
(\gr\mathcal{O}_{\mathcal{M}})_1$. If
$\psi=(\psi_{\red},\psi^*):(M,\mathcal{O}_{\mathcal{M}})\to (N,\mathcal{O}_{\mathcal{N}})$
is a morphism of supermanifolds, then $\gr(\psi)=(\psi_{\red},\gr(\psi^*))$, where
$\gr(\psi^*):\gr \mathcal{O}_{\mathcal{N}} \to \gr \mathcal{O}_{\mathcal{M}}$ is defined by
$$
\gr(\psi^*)(f+\mathcal{J}_{\mathcal{N}}^p): = \psi^*(f)+\mathcal{J}_{\mathcal{M}}^p
\,\,\text{for}\,\, f\in (\mathcal{J}_{\mathcal{N}})^{p-1}.
$$
Recall that by definition every morphism $\psi$ of supermanifolds is
even and as consequence sends $\mathcal{J}_{\mathcal{N}}^p$ into
$\mathcal{J}_{\mathcal{M}}^p$.

\medskip

\de
The supermanifold $\gr(\mathcal{M})$ is called the {\it
retract} of $\mathcal{M}$.
\medskip

For classification of supermanifolds we use
the following corollary of the well-known Green Theorem (see \cite{Green}, \cite{Bunegina_Oni} or \cite{DWitt} for more details):

\begin{thm}\label{Green corollary}[Green] {\it Let ${\mathcal{N}} =
(\mathcal{N}_{0}, \bigwedge \mathcal{E})$ be a split supermanfold of
dimension $n|m$, where $m\le 3$. The classes of isomorphic
supermanifolds $\mathcal{M}$ such that  $\gr \mathcal{M} =
{\mathcal{N}}$ are in bijection with orbits of the action
${\mathrm{Int}}$ of the group ${\mathrm {Aut}}\,{\mathbf E}$ on
$H^1(\mathcal{M}_0,{(\mathcal T}_{{\mathcal{N}}})_2)$.}

\end{thm}
\medskip

\noindent{\bf Remark.} This theorem allows to classify
supermanifolds $\mathcal{M}$ such that $\gr \mathcal{M} =
{\mathcal{N}}$ is fixed up
to isomorphisms that induce the identity morphism on $\gr
\mathcal{M}$.

\medskip

\section{Supermanifolds associated with $\mathbb{CP}^1$}

In what follows we will consider supermanifolds with the underlying space
${\mathbb{CP}}^1$.

\subsection{Supermanifolds with underlying space $\mathbb{CP}^1$}

 Let $\mathcal{M}$ be a
super\-mani\-fold of dimension $1|m$. Denote by $U_0$ and $U_1$ the
standard charts on ${\mathbb{CP}}^1$ with coordinates  $x$ and
$y=\frac1x$, respectively. By the Birkhoff-Grothendieck Theorem we can cover
$\gr\mathcal{M}$ by two charts
$$
(U_0,\mathcal{O}_{\gr
\mathcal{M}}|_{U_0})\,\,\text{ and}\,\,\, (U_1,\mathcal{O}_{\gr
\mathcal{M}}|_{U_1})
$$
with local coordinates
$x,\,\xi_1,\ldots,\xi_m$ and $y,\,\eta_1,\ldots,\eta_m$,
respectively, such that in $U_0\cap U_1$ we have
$$
y = x^{-1},\quad \eta_i = x^{-k_i}\xi_i,\; i = 1,\ldots,m,
$$
where $k_i$ are integers. Note that a
permutation of $k_i$ induces the automorphism of $\gr\mathcal{M}$.
  We will identify
$\gr\mathcal{M}$ with the set $(k_1,\ldots, k_m)$, so we will say that a
 supermanifold has the retract $(k_1,\ldots, k_m)$. In this paper
 we study the case: $m=3$ and $k_1=k_2=k_3=:k$. From now on   we use the
  notation $\mathcal{T} = \bigoplus \mathcal{T}_p$ for the tangent sheaf of $\gr\mathcal{M}$.

\subsection{A basis of $H^1({\mathbb{CP}}^1,{\mathcal T}_2).$}

Assume
that $m = 3$ and that $\mathcal{M} = (k_1,k_2,k_3)$ is a split super\-mani\-fold
with the underlying space $\mathcal{M}_0=\mathbb{CP}^1$. Let ${\mathcal T}$ be
its tangent sheaf. In \cite{Bunegina_Oni} the following decomposition
\begin{equation}\label{s-invar decomp}
{\mathcal T}_2=\sum_{i<j}{\mathcal T}_2^{ij}
\end{equation}
was obtained. The sheaf ${\mathcal T}_2^{ij}$ is a locally free
sheaf of rank $2$; its basis sections over
$(U_0,\mathcal{O}_{\mathcal{M}}|_{U_0})$ are:
\begin{equation}\label{loc section of T^ij}
\xi_i\xi_j\frac{\partial}{\partial x}, \ \
\xi_i\xi_j\xi_l\frac{\partial}{\partial \xi_l};
\end{equation}
where $l\ne i,j$. In $U_0\cap U_1$ we have
\begin{equation}\label{sootnoshen v peresech}
\begin{aligned}
\xi_i\xi_j\frac{\partial}{\partial x}
&=-y^{2-k_i-k_j}\eta_i\eta_j\frac{\partial}{\partial y} -
k_l y^{1-k_i-k_j}\eta_i\eta_j\eta_l\frac{\partial}{\partial \eta_l},\\
\xi_i\xi_j\xi_l\frac{\partial}{\partial \xi_l}
&= y^{-k_i-k_j}\eta_i\eta_j\eta_l\frac{\partial}{\partial \eta_l}.
\end{aligned}
\end{equation}

The following theorem was proven in \cite{Vish}. For completeness we reproduce it here.

\begin{thm}\label{basis of H^1} Assume that $i<j$ and $l\ne i,j$.

$\bf1.$  For  $k_i+k_j>3$ the basis
of $H^1({\mathbb{CP}}^1,{\mathcal T}_2^{ij})$ is:

\begin{equation}\label{cocycles generating H^1}
\begin{aligned}
&x^{-n}\xi_i\xi_j\frac{\partial}{\partial x},\ \ n=1,{\ldots},k_i+k_j-3,\\
&x^{-n}\xi_i\xi_j\xi_l\frac{\partial}{\partial \xi_l},\ \
n=1,{\ldots},k_i+k_j-1;
\end{aligned}
\end{equation}

$\bf2.$ for  $k_i+k_j=3$ the basis
of $H^1({\mathbb{CP}}^1,{\mathcal T}_2^{ij})$ is:
$$
x^{-1}\xi_i\xi_j\xi_l\frac{\partial}{\partial \xi_l},\quad
x^{-2}\xi_i\xi_j\xi_l\frac{\partial}{\partial \xi_l};
$$

$\bf3.$ for $k_i+k_j=2$ and $k_l=0$ the basis
of $H^1({\mathbb{CP}}^1,{\mathcal T}_2^{ij})$ is:
$$
x^{-1}\xi_i\xi_j\xi_l\frac{\partial}{\partial \xi_l};
$$

$\bf4.$ if $k_i+k_j=2$ and $k_l\ne 0$ or $k_i+k_j<2$, we have
$H^1({\mathbb{CP}}^1,{\mathcal T}_2^{ij})=\{0\}$.

\end{thm}
\medskip

\noindent{\it Proof.}  We use the \v{C}ech cochain complex of the cover
${\mathfrak U} =\{U_0,U_1\}$ as in $3.1$, hence, $1$-cocycle with values in the
sheaf ${\mathcal T}_2^{ij}$ is a section $v$ of ${\mathcal
T}_2^{ij}$ over $U_0\cap U_1$. We are looking for {\it basis
cocycles}, i.e. cocycles such that their cohomology classes form a
basis of
 $H^1({\mathfrak U},{\mathcal
T}_2^{ij})\simeq H^1({\mathbb{CP}}^1,{\mathcal T}_2^{ij})$. Note
that if $v\in Z^1({\mathfrak U},{\mathcal T}_2^{ij})$ is holomorphic
in $U_0$ or $U_1$ then the cohomology class of $v$ is equal to $0$.
Obviously, any $v\in Z^1({\mathfrak U},{\mathcal T}_2^{ij})$ is a
linear combination of vector fields (\ref{loc section of T^ij}) with
holomorphic in $U_0\cap U_1$ coefficients. Further, we expand these
coefficients in a Laurent series in $x$ and drop the summands
$x^n,\; n\ge 0$, because they are holomorphic in $U_0$. We see that
$v$ can be replaced by
\begin{equation}\label{cocycle chomological to v}
v = \sum_{n=1}^{\infty} a^n_{ij}x^{-n}\xi_i\xi_j\frac{\partial}
{\partial x} + \sum_{n=1}^{\infty}
b^n_{ij}x^{-n}\xi_i\xi_j\xi_l\frac{\partial} {\partial \xi_l},
\end{equation}
where $a^n_{ij}, b^n_{ij}\in{\mathbb C}$. Using (\ref{sootnoshen v
peresech}), we see that the summands corresponding to $n\ge
k_i+k_j-1$ in the first sum of (\ref{cocycle chomological to v}) and
the summands corresponding to $n\ge k_i+k_j$ in the second sum of
(\ref{cocycle chomological to v}) are holomorphic in $U_1$. Further,
it follows from (\ref{sootnoshen v peresech}) that
$$
x^{2-k_i-k_j}\xi_i\xi_j\frac{\partial}{\partial x}\sim
-k_lx^{1-k_i-k_j}\xi_i\xi_j\xi_l\frac{\partial}{\partial \xi_l}.
$$
Hence the cohomology classes of the following cocycles
$$
\begin{aligned}
&x^{-n}\xi_i\xi_j\frac{\partial}{\partial x},\ \ n=1,{\ldots},k_i+k_j-3,\\
&x^{-n}\xi_i\xi_j\xi_l\frac{\partial}{\partial \xi_l},\ \
n=1,{\ldots},k_i+k_j-1,
\end{aligned}
$$
generate $H^1({\mathbb{CP}}^1,{\mathcal T}_2^{ij})$. If we examine
linear combination of these cocycles which are
cohomological trivial, we get the result.$\Box$

\medskip

\noindent{\bf Remark.}
Note that a similar method can be used for computation of a basis of
$H^1({\mathbb{CP}}^1,{\mathcal T}_q)$ for any odd dimension $m$ and any $q$.

\medskip

In the case $k_1=k_2=k_3=k$, from Theorem \ref{basis of H^1}, it follows:

\begin{cor}\label{basis of H^1 (k,k,k)}  Assume that $i<j$ and $l\ne i,j$.

 $\bf1.$ For  $k\geq 2$ the basis
of $H^1({\mathbb{CP}}^1,{\mathcal T}_2^{ij})$ is

\begin{equation}\label{cocycles generating H^1(k,k,k)}
\begin{aligned}
&x^{-n}\xi_i\xi_j\frac{\partial}{\partial x},\ \ n=1,{\ldots},2k-3,\\
&x^{-n}\xi_i\xi_j\xi_l\frac{\partial}{\partial \xi_l},\ \
n=1,{\ldots},2k - 1.
\end{aligned}
\end{equation}

$\bf2.$ If $k<2$, we have $H^1({\mathbb{CP}}^1,{\mathcal T}_2)=\{0\}$.

\end{cor}

\subsection{The group $\operatorname{Aut}\bold E$}

This section is devoted to the calculation of the group of automorphisms
$\operatorname{Aut}\bold E$ of the vector bundle $\bold E$ in the case $(k,k,k)$.
 Here $\bold E$ is the vector bundle corresponding to the split supermanifold $(k,k,k)$.

Let $(\xi_i)$ be a local basis of $\bold E$ over $U_0$ and $A$ be an
automorphism of $\bold E$. Assume that $A(\xi_j)=\sum
a_{ij}(x)\xi_i$. In $U_1$ we have
$$
A(\eta_j) = A(y^{k_j}\xi_j) = \sum y^{k_j- k_i}a_{ij}(y^{-1})\eta_i.
$$
Therefore, $a_{ij}(x)$ is a polynomial in $x$ of degree $\leq
k_j-k_i$, if
 $k_j-k_i\ge 0$ and is $0$, if $k_j-k_i<0$.
We denote by $b_{ij}$ the entries of the matrix $B=A^{-1}$. The
entries are also polynomials in $x$ of degree $\leq k_j-k_i$. We
 need the following formulas:
\begin{equation}\label{A(vector field)}
\begin{array}{c}
A(\xi_1\xi_2\xi_3\frac{\partial}{\partial \xi_k})A^{-1} =
\det(A)\sum\limits_s
b_{ks} \xi_1\xi_2\xi_3 \frac{\partial}{\partial \xi_s};\\
\rule{0pt}{6mm}A(\xi_i\xi_j\frac{\partial}{\partial x})A^{-1} =
\det(A) \sum\limits_{k<s} (-1)^{l+r}b_{lr}
\xi_k\xi_s\frac{\partial}{\partial x} + \\
+\det(A)\sum\limits_s b'_{ls}\xi_i\xi_j\xi_l\frac{\partial}{\partial
\xi_s},
\end{array}
\end{equation}
where $i < j$, $l\ne i,j$ and $r\ne k,s$. Here we use the notation $b'_{ls} =
\frac{\partial}{\partial x}(b_{ls})$.
By (\ref{A(vector field)}), in the case $k_1=k_2=k_3=k$, we have:

\begin{prop}\label{Aut (k,k,k)} Assume that $k_1=k_2=k_3=k$.

${\bf 1.}$ We have
$$
\operatorname{Aut}\bold E \simeq \GL_3(\mathbb{C}).
$$
 In other words
$$
\operatorname{Aut}{\mathbf{E}} =\{ (a_{ij})\,|\, a_{ij}\in  \mathbb{C},\, \det(a_{ij})\ne 0 \}.
$$

${\bf 2.}$ The action of $\operatorname{Aut}\bold E $ on $\mathcal
{T}_2$ is given in $U_0$ by the following formulas:

\begin{equation}\label{A(vector field){k,k,k}}
\begin{array}{c}
A(\xi_1\xi_2\xi_3\frac{\partial}{\partial \xi_k})A^{-1} =
\det(A) \xi_1\xi_2\xi_3\sum\limits_s
b_{ks} \frac{\partial}{\partial \xi_s};\\

\rule{0pt}{6mm}A(\xi_i\xi_j\frac{\partial}{\partial x})A^{-1} =
\det(A) \sum\limits_{k<s} (-1)^{l+r}b_{lr}
\xi_k\xi_s\frac{\partial}{\partial x},
\end{array}
\end{equation}
where $i < j$, $l\ne i,j$ and $r\ne k,s$. Here $B = (b_{ij}) = A^{-1}$
\end{prop}

\section{Classification of supermanifolds with retract $(k,k,k)$}

In this section we give a classification up to isomorphism of complex-analytic supermanifolds
with underlying space ${\mathbb{CP}}^1$ and with retract $(k,k,k)$ using Theorem
\ref{Green corollary}. In previous section we have calculated the vector space
$H^1({\mathbb{CP}}^1,{\mathcal T}_2)$, the group $\operatorname{Aut}\bold E$ and
the action ${\operatorname{Int}}$ of $\operatorname{Aut}\mathbf{E}$  on
$H^1({\mathbb{CP}}^1,{\mathcal T}_2)$, see Theorem \ref{basis of H^1 (k,k,k)} and Proposition \ref{Aut (k,k,k)}.
Our objective in this section is to calculate the orbit space  corresponding to the action ${\operatorname{Int}}$:
\begin{equation}
\label{orbit space}
H^1({\mathbb{CP}}^1,{\mathcal T}_2)/ \operatorname{Aut}\mathbf{E}.
\end{equation}
By Theorem  \ref{Green corollary} classes of isomorphic supermanifolds are in one-to-one correspondence with points of the set (\ref{orbit space}).

Let us fix the following basis of $H^1({\mathbb{CP}}^1,{\mathcal T}_2)$:
\begin{equation}\label{basis1}
\begin{array}{lll}
v_{11} = x^{-1}\xi_2\xi_3\frac{\partial}{\partial x},&v_{12} = -x^{-1}\xi_1\xi_3\frac{\partial}{\partial x}, &
v_{13} = x^{-1}\xi_1\xi_2\frac{\partial}{\partial x},\\
\cdots & \cdots&  \cdots\\
v_{p1} = x^{-p}\xi_2\xi_3\frac{\partial}{\partial x},&v_{p2} = -x^{-p}\xi_1\xi_3\frac{\partial}{\partial x}, &
v_{p3} = x^{-p}\xi_1\xi_2\frac{\partial}{\partial x},\\
\end{array}
\end{equation}
\begin{equation}\label{basis2}
\begin{array}{lll}
v_{p+1,1} = x^{-1}\xi_1\xi_2\xi_3\frac{\partial}{\partial \xi_1},&\cdots & v_{p+1,3} = x^{-1}\xi_1\xi_2\xi_3\frac{\partial}{\partial \xi_3},\\

\cdots & \cdots&  \cdots\\

v_{q1} = x^{-a}\xi_1\xi_2\xi_3\frac{\partial}{\partial \xi_1},&\cdots & v_{q3} = x^{-a}\xi_1\xi_2\xi_3\frac{\partial}{\partial \xi_3},
\end{array}
\end{equation}
where $p=2k-3$, $a= 2k-1$ and $q= p+a= 4k-4$. (Compere with Theorem \ref{basis of H^1 (k,k,k)}.)  Let us take
$A\in \operatorname{Aut}\mathbf{E} \simeq\GL_3({\mathbb{C}})$, see Proposition \ref{Aut (k,k,k)}.
We get that in the basis (\ref{basis1}) - (\ref{basis2})
 the
automorphism ${\operatorname{Int}}\, A$ is given by
$$
{\operatorname{Int}}\, A(v_{is}) = \frac{1}{ \det B}\sum_jb_{sj}v_{ij}.
$$
 Note that for any matrix $C\in
\GL_3(\mathbb C)$ there exists a matrix $B$ such that
$$
C=\frac1{\det\,B}B.
$$ Indeed, we can put
$B=\frac{1}{\sqrt{\det\,C}}C$. We summarize these observations in the following proposition:

\begin{prop}\label{lem_main} Assume that $k_1=k_2=k_3=k.$ Then
$$
H^1({\mathbb{CP}}^1,{\mathcal T}_2) \simeq \operatorname{Mat}_{3\times (4k-4)}(\mathbb{C})
$$
 and the action ${\operatorname{Int}}$ of $\operatorname{Aut}\mathbf{E}$  on
 $H^1({\mathbb{CP}}^1,{\mathcal T}_2)$  is equivalent to the standard
  action of $\GL_3(\mathbb{C})$ on $ \operatorname{Mat}_{3\times (4k-4)}(\mathbb{C})$. More precisely,
   ${\operatorname{Int}}$ is equivalent to the following action:
\begin{equation}
\label{action}
 D\longmapsto (W \longmapsto DW),
\end{equation}
where $D\in \GL_3(\mathbb{C})$, $W\in \operatorname{Mat}_{3\times (4k-4)}(\mathbb{C})$ and $DW$ is the usual matrix multiplication.
\end{prop}

Now we prove our main result.

\begin{thm}\label{teor_main} Let $k \geq 2$.  Complex-analytic supermanifolds with underlying space ${\mathbb{CP}}^1$ and retract $(k,k,k)$ are
in one-to-one correspondence with points of the following set:
$$
\bigcup\limits_{r=0}^3\mathbf{Gr}_{4k-4,r},
$$
where $\mathbf{Gr}_{4k-4,r}$ is  the Grassmannian of type $(4k-4,r)$, i.e. it is the set of all $r$-dimensional subspaces in $\mathbb{C}^{4k-4}$.

In the case $k<2$ all supermanifolds with underlying space ${\mathbb{CP}}^1$ and retract $(k,k,k)$ are split and isomorphic to their retract $(k,k,k)$.

\end{thm}

\medskip

\noindent{\it Proof.} Assume that $k \geq 2$. Clearly, the action (\ref{action}) preserves
the rank $r$ of matrices from $\operatorname{Mat}_{3\times
(4k-4)}(\mathbb{C})$ and $r\leq 3$. Therefore, matrices with
different rank belong to different orbits of this action.
Furthermore, let us fix $r\in \{ 0,1,2,3\}$. Denote by
$\operatorname{Mat}^r_{3\times (4k-4)}(\mathbb{C})$ all matrices
with rank $r$.  Clearly, we have
$$
\operatorname{Mat}_{3\times (4k-4)}(\mathbb{C}) = \bigcup\limits_{r=0}^3 \operatorname{Mat}^r_{3\times (4k-4)}(\mathbb{C}).
$$

A matrix $W \in \operatorname{Mat}^r_{3\times (4k-4)}(\mathbb{C})$
defines the $r$-dimensional subspace $V_W$ in $\mathbb{C}^{4k-4}$.  This subspace is the
linear combination of lines of $W$. (We consider lines of a matrix $X\in \operatorname{Mat}_{3\times (4k-4)}(\mathbb{C})$
as vectors in  $\mathbb{C}^{4k-4}$.) Therefore, we have defined the map $F_r$:
$$
W\longmapsto F_r(W)=V_W\in \mathbf{Gr}_{4k-4,r}.
$$

The map $F_r$ is surjective. Indeed, in any $r$-dimensional subspace $V\in \mathbf{Gr}_{4k-4,r}$, where $r\leq 3$,
we can take $3$ vectors generating $V$ and form the matrix $W_V\in \operatorname{Mat}^r_{3\times (4k-4)}(\mathbb{C})$.
In this case the matrix $W_V$ is of rank $r$ and $F_r(W_V)=V$. Clearly, $F_r(W)= F_r(DW)$, where $D\in  \GL_3(\mathbb{C})$.
 Conversely, if $W$ and $W'\in F_r^{-1}(V_W)$, then there exists a matrix $D\in \GL_3(\mathbb{C})$ such that $DW= W'$.
  It follows that orbits of $\GL_3(\mathbb{C})$ on  $\operatorname{Mat}^r_{3\times (4k-4)}(\mathbb{C})$
  are in one to one correspondence with points of $\mathbf{Gr}_{4k-4,r}$. Therefore, orbits of $\GL_3(\mathbb{C})$ on
$$
\operatorname{Mat}_{3\times (4k-4)}(\mathbb{C})= \bigcup_{r=0}^3 \operatorname{Mat}^r_{3\times (4k-4)}(\mathbb{C})
$$
are in one-to-one correspondence with points of $\bigcup\limits_{r=0}^3\mathbf{Gr}_{4k-4,r}$. The proof is\linebreak
 complete.$\Box$

\medskip

\section{Appendix. Classification of supermanifolds\\
 with underlying space $\mathbb{CP}^1$ of odd dimension $2$. }

In this Section we give a classification up to isomorphism of complex-analytic supermanifolds in the case $m=2$ and $\gr \mathcal{M} =
(k_1,k_2)$, where $k_1$, $k_2$ are any integers. As far as we know the classification in this case
does not appear in the literature, but it is known for specialists.

 Let us compute  the 1-cohomology with values in the
tangent sheaf $\mathcal{T}_2$. The sheaf $\mathcal{T}_2$ is a
locally free sheaf of rank $1$.  Its basis section over $(U_0,
\mathcal{O}_{\mathcal{M}}|_{U_0})$ is
$\xi_1\xi_2\frac{\partial}{\partial x}.$ The transition functions in
$U_0 \cap U_1$ are given by the following formula:
$$
\xi_1\xi_2\frac{\partial}{\partial x} = -y^{2-k_1-k_2} \eta_1\eta_2\frac{\partial}{\partial y}.
$$
Therefore, a basis of $H^1(\mathbb{CP}^1,\mathcal{T}_2)$ is
$$
x^{-n}\xi_1\xi_2\frac{\partial}{\partial x}, \,\,n = 1,\cdots, k_1+k_2-3.
$$

Let $(\xi_i)$ be a local basis of $\bold E$ over $U_0$ and $A$ be an
automorphism of $\bold E$. As in the case $m=3$, we have that $a_{ij}(x)$ is a polynomial in $x$ of degree $\leq
k_j-k_i$, if
 $k_j-k_i\ge 0$ and is $0$, if $k_j-k_i<0$. We need the following
 formulas:
 $$
 A(x^{-n}\xi_1\xi_2\frac{\partial}{\partial x})A^{-1} = (\det A) x^{-n}\xi_1\xi_2\frac{\partial}{\partial x}.
 $$
Denote
$$
v_n= x^{-n}\xi_1\xi_2\frac{\partial}{\partial x}, \,\,n = 1,\cdots, k_1+k_2-3.
$$
We see that the action ${\mathrm{Int}}$ is equivalent to the action of $\mathbb{C}^*$ on $\mathbb{C}^{k_1+k_2-3}$, therefore, the quotient
space is $\mathbb{CP}^{k_1+k_2-4}\cup \{\mathrm{pt}\}$ for $k_1+k_2 \geq 4$ and $\{\mathrm{pt}\}$ for $k_1+k_2<4$. We have proven the following theorem:

\begin{thm}\label{teor_main m=2} Assume that $k_1+k_2 \geq 4$.  Complex-analytic supermanifolds
with underlying space ${\mathbb{CP}}^1$ and retract $(k_1,k_2)$ are
in one-to-one correspondence with points of
$$
\mathbb{CP}^{k_1+k_2-4}\cup \{\mathrm{pt}\}.
$$

In the case $k_1+k_2<4$ all supermanifolds with underlying space ${\mathbb{CP}}^1$ and
retract $(k_1,k_2)$ are split and isomorphic to their retract  $(k_1,k_2)$.

\end{thm}

\subsection*{Acknowledgment}
The author is grateful to A.\,L.\,Onishchik for useful discussions.

\noindent{\it Elizaveta Vishnyakova}

\noindent University of Luxembourg

 \noindent {\emph{E-mail address:}
\verb"VishnyakovaE@googlemail.com"}


\begin{thebibliography}{1}
\bibitem[B]{Ber} {\it  Berezin F.A.} Introduction to superanalysis. Edited and with a foreword by A. A. Kirillov.
With an appendix by V. I. Ogievetsky. Mathematical Physics and
Applied Mathematics, 9. D. Reidel Publishing Co., Dordrecht, 1987.

\bibitem[BL]{BL} {\it Berezin F.A., Leites D.A.} Supermanifolds. Soviet
Math. Dokl. 16, 1975, 1218-1222.


\bibitem[BO]{Bunegina_Oni} {\it Bunegina V.A. and Onishchik A.L.} Homogeneous supermanifolds associated with
the complex projective line. Algebraic geometry, 1. J. Math. Sci. 82
(1996), no. 4, 3503-3527.

\bibitem[DW]{DWitt} {\it Donagi R. and Witten E.}   Supermoduli space is not projected. arXiv:1304.7798, 2013.



\bibitem[Gr]{Green} {\it Green P.} On holomorphic graded manifolds. Proc. Amer. Math. Soc. 85 (1982), no. 4, 587-590.

\bibitem[L]{ley} {\it Leites D.A.} Introduction to the theory of supermanifolds.
Russian Math. Surveys 35 (1980), 1-64.

\bibitem[V]{Vish}  {\it Vishnyakova E. G.} Even-homogeneous supermanifolds on the complex projective
line.    Differential Geometry and its Applications 31 (2013) 698-706.


\end{thebibliography}
\end{document}